\documentclass[11pt,a4paper]{article}
\usepackage{amsmath,mathtools}
\usepackage{amssymb}
\newtheorem{Thm}{Theorem}[section]
\newtheorem{Lem}[Thm]{Lemma}

\newtheorem{Prop}[Thm]{Proposition}

\newtheorem{Def}[Thm]{Definition}
\newtheorem{Ex}[Thm]{Example}
\numberwithin{equation}{section}

\newcommand{\Z}{\mathbb{Z}}

\newcommand{\N}{\mathbf{N}}
\newcommand{\R}{\mathbb{R}}
\newcommand{\C}{\mathbb{C}}

\newcommand{\T}{\mathbb{T}}

\newcommand{\bpr}{ \textit{Proof}: }

\newcommand{\eprsk}{~$\blacksquare$\medskip}

\newcommand{\SL}{\operatorname{SL}}

\providecommand{\AMS}{$\mathcal{A}$\kern-.1667em%
\lower.25em\hbox{$\mathcal{M}$}\kern-.125em$\mathcal{S}$}
\begin{document}
\date{May 6, 2024}
\title{Amenable actions of real and $p$-adic algebraic groups}

\author{Alain VALETTE}

\maketitle


\begin{abstract} Let $K$ be a locally compact field of characteristic 0. Let $G$ be a linear algebraic group defined over $K$, acting algebraically on an algebraic variety $V$. We prove that the action of $G(K)$ (the group of $K$-rational points of $G$) on $V(K)$ is topologically amenable, if and only if all points stabilizers in $G(K)$ are solvable-by-compact. This follows by combining a result by Borel-Serre \cite{BoSe} with the following fact: let $G$ be a second countable locally compact group acting continuously on a second countable locally compact space $Y$. If the action $G\curvearrowright Y$ is smooth (i.e. the Borel structure on $G\backslash Y$ is countably separated), then topological amenability of $G\curvearrowright Y$ is equivalent to amenability of all point stabilizers in $G$. 
\end{abstract}
\section{Introduction}

Amenability of group actions is a far-reaching generalization of classical amenability for locally compact groups. It was proposed by R.J. Zimmer (Def. 4.3.1 in \cite{Zim}) in the measurable setting, and in the topological setting by J. Renault (measurewise amenability, see Def. 2.3.6 in \cite{Ren}) and C. Anantharaman-Delaroche (topological amenability, see Def. 2.1 in \cite{Anan}). We recall the relevant definitions.

\begin{Def}\label{MeasureAmen}\footnote{This is not Zimmer's original definition, but it is equivalent to it by Theorem A in \cite{AEG} .}Let $G$ be a locally compact group acting measurably on a measure space $(X,\mu)$, where $\mu$ is a quasi-invariant measure. We say that the action $G\curvearrowright X$ is {\bf amenable in the sense of Zimmer} if there exists a $G$-equivariant conditional expectation $L^\infty(X\times G,\mu\times\lambda_G)\rightarrow L^\infty(X,\mu)$, where $\lambda_G$ denotes Haar measure on $G$.
\end{Def}

\begin{Def}\label{TopAmen} Let $G$ be a locally compact group acting continuously on a locally compact space $Y$. 
\begin{enumerate}
\item[1)] The action $G\curvearrowright Y$ is {\bf measurewise amenable} if the action of $G$ on $(Y,\mu)$ is amenable in the sense of Zimmer for every quasi-invariant measure $\mu$ on $Y$.
\item[2)] We denote by $Prob(G)$ the set of probability measures on $G$, equipped with the weak $^*$-topology. We say that the action $G\curvearrowright Y$ is {\bf topologically amenable} if there exists a net $(m^i)_{i\in I}$ of continuous functions $m^i:Y\rightarrow Prob(G): y\mapsto m^i_y$ such that $\lim_{i\rightarrow\infty}\|g_*m^i_y - m^i_gy\|=0$ uniformly on compact subsets of $Y\times G$ (where the norm corresponds to the total variation distance).

\end{enumerate}
\end{Def}

It follows from any of the three definitions above, that a locally compact group is amenable if and only if it acts amenably on a one-point space.

It was proved by Anantharaman-Delaroche and Renault that topological amenability implies measurewise amenability (Proposition 3.3.5 in \cite{ADR}), and that they are equivalent for actions of countable discrete groups. The equivalence for general group actions was a long-standing open question that was settled only recently through combined work by Buss-Echterhoff-Willett \cite{BEW} and Bearden-Crann \cite{BeCr}: {\it for a locally compact second countable group $G$ acting countinuously on a locally compact second countable space $Y$, the action $G\curvearrowright Y$ is topologically amenable if and only if it is measurewise amenable} (see Corollary 3.29 in \cite{BEW}).

The goal of this Note is to prove:

\begin{Thm}\label{main} Let $K$ be a locally compact field of characteristic 0. Let $G$ be a linear algebraic group defined over $K$, acting algebraically on an algebraic variety $V$. The action of $G(K)$ (the group of $K$-rational points of $G$) on $V(K)$ is topologically amenable, if and only if all points stabilizers in $G(K)$ are solvable-by-compact.
\end{Thm}

Recall that, for a locally compact second countable group $G$ acting continuously on a locally compact second countable space $Y$, the action of $G$ on $Y$ is {\it smooth} if the quotient Borel structure on $G\backslash Y$ is countably separated. In section 2 we prove that, for smooth actions, topological amenability is equivalent to amenability of stabilizers (see Proposition \ref{smooth}).
As explained in section 3, Theorem \ref{main} follows then easily by combining this result with a deep result by Borel-Serre \cite{BoSe} asserting that the action $G(K)\curvearrowright V(K)$ is always smooth.


We end this Introduction with the example that motivated the present study.

\begin{Ex} For $n\geq 0$, let $\rho_n$ denote the $(n+1)$-dimensional irreducible representation of $SL_2(K)$ on the space $P_n(K)$ of homogeneous polynomials of degree $n$ with coefficients in $K$ in the 2 variables $X,Y$. So, for $P(X,Y)\in P_n(K)$ and $\left(\begin{array}{cc}a & b \\c & d\end{array}\right)\in SL_2(K)$, we define
$$(\rho_n(\left(\begin{array}{cc}a & b \\c & d\end{array}\right))P)(X,Y) = P(aX+cY,bX+dY).$$
Going to the projective space of $P_n(K)$, we get an action of $SL_2(K)$ on the projective space $\mathbb{P}^n(K)$. In Proposition 3.5 of \cite{JoVa}, the authors gave a direct but {\it ad hoc} proof of the fact that, for $K=\R,\C$, the action $SL_2(K)\curvearrowright \mathbb{P}^n(K)$ is topologically amenable. Now, for any local field $K$ of characteristic 0, stabilizers of points in $\mathbb{P}^n(K)$ are proper algebraic subgroups of dimension at most 2, so they are solvable. By Theorem \ref{main}, the action $SL_2(K)\curvearrowright \mathbb{P}^n(K)$ is topologically amenable.

\end{Ex}


\section{Amenable actions vs. amenable equivalence relations}

The next lemma is probably well-known.
\begin{Lem}\label{amenstab} Let the locally compact group $G$ act continuously on a locally compact space $Y$. If the action of $G$ on $Y$ is topologically amenable, then for every $y\in Y$ the stabilizer $G_y$ is amenable.
\end{Lem}

\bpr As the action $G\curvearrowright Y$ is topologically amenable, by Proposition 2.5 in \cite{Anan} there exists a net $(h_i)_{i\in I}$ of continuous, compactly supported, positive type functions\footnote{Recall that a function $h$ on $Y\times G$ is positive type if, for every $n\geq 1, y\in Y, t_1,...,t_n\in G,\lambda_1,...,\lambda_n\in\C$, we have $\sum_{i,j}\overline{\lambda_i}\lambda_j h(t_i^{-1}y,t_i^{-1}t_j)\geq 0$.} on $Y\times G$ that converges to 1 uniformly on compact subsets on $Y\times G$. For $t\in G_y$, set $g_i(t)=h_i(y,t)$. Then $(g_i)_{i\in I}$ is a net of continuous, compactly supported, positive type functions on $G_y$ that converges to 1 on compact subsets of $G_y$: the existence of such a net characterizes amenability of $G_y$. \eprsk

\bigskip
We must now discuss amenable equivalence relations as introduced in \cite{CFW}. So let $X$ be a standard Borel space and $\mathcal{R}\subset X\times X$ be a Borel equivalence relation. A transverse measure $\Lambda$ for $\mathcal{R}$ is defined as in Definition 1 in Chapter II in \cite{Con}; our aim is to define the amenability of the pair $(\mathcal{R},\mbox{class of}\,\Lambda)$. Denote by $p_1:X\times X\rightarrow X:(x,y)\mapsto x$ the first projection map. We choose an auxiliary {\it transverse function} $\nu$ for $\mathcal{R}$, i.e. a map $x\mapsto \nu^x$ from $X$ to the space of positive measures on $\mathcal{R}$, such that:
\begin{itemize}
\item For every $x\in X$, the measure $\nu^x$ is non-zero and supported on $[x]_\mathcal{R}=:p_1^{-1}(x)\cap\mathcal{R}$. Moreover $\nu$ is invariant in the sense that $\gamma\nu^x=\nu^y$ for every $\gamma=(y,x)\in\mathcal{R}$.
\item For every Borel set $A$ in $\mathcal{R}$, the map $X\rightarrow [0,+\infty]: x\mapsto \nu^x(A)$ is measurable.
\item $\nu$ is proper in the sense that $\mathcal{R}$ is a countable union of Borel sets $(A_n)_{n\in\N}$ such that, for every $n\in\N$, the function $x\mapsto \nu^x(A_n)$ is bounded.
\end{itemize}
Theorem 3 in Chapter II of \cite{Con} then provides a bijection between transverse measures on $\mathcal{R}$ and ordinary measures $\mu$ on $X$ such that the measure $m=:\int_X\nu^x\,d\mu(x)$ on $\mathcal{R}$ is invariant under the flip on $\mathcal{R}$.

\begin{Def}\label{amenrel}[See p. 446 in \cite{CFW}] The pair $(\mathcal{R},\mbox{class of}\,\Lambda)$ is amenable if, for almost every $x\in X$ there exists a state $p_x$ on $L^\infty([x]_\mathcal{R},\nu^x)$ such that:
\begin{enumerate}
\item The family $(p_x)_{x\in X}$ is invariant in the sense that $\gamma p_x=p_y$ for every $\gamma=(y,x)\in\mathcal{R}$.
\item With $m=\int_X\nu^x\,d\mu(x)$ as above, for every $f\in L^\infty(\mathcal{R},m)$ the function $x\mapsto p_x(f)$ is $\mu$-measurable.
\end{enumerate}
\end{Def}

It is proved in lemma 15 of \cite{CFW} that this definition is in fact independent of the choice of the transverse function $\nu$. 

With $X,\mathcal{R}$ as above, we say that the equivalence relation $\mathcal{R}$ is {\it smooth} if the quotient Borel structure on $\mathcal{R}\backslash X$ is countably separated (see Definition 2.1.9 in \cite{Zim}). Extending the first part of Definition \ref{TopAmen}, we say that $\mathcal{R}$ is {\it measurewise amenable} if $(\mathcal{R},\mbox{class of}\,\Lambda)$ is amenable for every transverse measure $\Lambda$. 

\begin{Lem}\label{smooth/amen} Any smooth Borel equivalence relation is measurewise amenable.
\end{Lem}

\bpr  Let $\mathcal{R}$ be a smooth equivalence relation on the standard Borel space $X$. \begin{itemize}
\item We denote by $\pi:X\rightarrow Y=:\mathcal{R}\backslash X$ the quotient map. By the discussion on page 47 of \cite{Con}, there is a bijection between transverse measures for $\mathcal{R}$ and ordinary measures on $Y$, so we choose a measure $\lambda$ on $Y$, which we may assume to be a probability measure. 
\item By the von Neumann selection theorem (Theorem A.9 in \cite{Zim}), there exists a co-null standard Borel set $Z\subset Y$ and a Borel section $s:Z\rightarrow Y$ such that $\pi\circ s=Id_{Z}$. Replacing $Y$ by $Z$ and $X$ by $\pi^{-1}(Z)$, we then define the Borel map $\iota: X\rightarrow \mathcal{R}: x\mapsto (x,s(\pi(x)))$ which is a section of $p_1|_{\mathcal{R}}$; observe that, for $\gamma=(y,x)\in\mathcal{R}$, we have $\gamma\iota(x)=\iota(y)$. We then define the transverse function $\nu$ as $\nu^x=:\delta_{\iota(x)}$ (the Dirac mass at $\iota(x)$, viewed as a probability measure on $\mathcal{R}$); we will use $\nu$ to check the conditions in Definition \ref{amenrel}. Observe that $\nu^{s(y)}=\delta_{\iota(s(y))}=\delta_{(s(y),s(y))}$ for $y\in Y$.
\item Set $\mu=:s_*\lambda$, form the measure $m=\int_X\nu^x\,d\mu(x)$ on $\mathcal{R}$. Since 
$$m=\int_Y\nu^{s(y)}\,d\lambda(y)=\int_Y \delta_{(s(y),s(y))}\,d\lambda(y),$$
the measure $m$ is supported on the diagonal, hence invariant under the flip on $\mathcal{R}$.
\item Finally, since $L^\infty ([x]_{\mathcal{R}},\nu^x)$ is one-dimensional, the state $p_x$ is tautological: it is given by $p_x(\varphi)=\varphi(\iota(x))$ for $\varphi\in L^\infty ([x]_{\mathcal{R}},\nu^x)$. The invariance of the family $(p_x)_{x\in  X}$ is then clear. 
Then for $f\in L^\infty(\mathcal{R},m)$, the function $x\mapsto p_x(f)=f(\iota(x))$ is clearly $\mu$-measurable. This concludes the proof.
\eprsk
\end{itemize}

Let $G$ be a locally compact second countable group acting continuously on a locally compact second countable space $Y$. The action of $G$ on $Y$ is {\it smooth} if the orbital equivalence relation $\mathcal{R}_G$ is smooth.
\begin{Prop}\label{smooth} Let $G$ act continuously on $Y$ as above, and assume that $G$ acts smoothly on $Y$. The following are then equivalent:
\begin{enumerate}
\item[(i)] The action $G\curvearrowright Y$ is topologically amenable.
\item[(ii)] All point stabilizers in $G$ are amenable.
\end{enumerate}
\end{Prop}

\bpr The implication $(i)\Rightarrow (ii)$ follows immediately from lemma \ref{amenstab}. For the converse implication $(ii)\Rightarrow (i)$: in view of the discussion following Definition \ref{TopAmen}, it is enough to prove that the action $G\curvearrowright Y$ is measurewise amenable. So fix a quasi-invariant measure $\mu$ on $Y$, we must prove that $G\curvearrowright (Y,\mu)$ is amenable in the sense of Zimmer. We appeal to a result by S. Adams, G. Elliott and T. Giordano (Theorem 5.1 in \cite{AEG}):  {\it Assume that the second countable locally compact group $G$ acts measurably on a standard measure space $(X, \nu)$ with $\nu$ a quasi-invariant probability measure. If stabilizers are amenable $\nu$-almost everywhere, and if the orbital equivalence relation $\mathcal{R}_G$ is amenable, then the action of $G$ on $(X,\nu)$ is amenable in the sense of Zimmer.}
In view of the Adams-Elliott-Giordano result, since stabilizers are assumed to be amenable, the result follows immediately from lemma \ref{smooth/amen}.
\eprsk


\section{Proof of Theorem \ref{main}}

We will appeal to a deep result of A. Borel and J.-P. Serre  (\cite{BoSe}, Proposition 4.9 and Corollaire 6.4; see also Theorem 3.1.3 in \cite{Zim}): all orbits of $G(K)$ in $V(K)$ are locally closed. By the Glimm-Effros theorem (see Theorem 2.1.14 in \cite{Zim}), having all orbits locally closed is equivalent to smoothness of the action, so we are in the assumptions of our Proposition \ref{smooth}, to the effect that topological amenability of $G(K)\curvearrowright V(K)$ is equivalent to amenability of stabilizers. Since solvable-by-compact groups are amenable, one implication in Theorem \ref{main} becomes clear. 

For the converse, assume that the action $G(K)\curvearrowright V(K)$ is amenable. Fix $y\in V(K)$, and denote by $H$ the stabilizer of $y$ in $G(K)$: by lemma \ref{amenstab}, the subgroup $H$ is amenable, and we must prove that $H$ is solvable-by-compact. Since $H$ is the group of $K$-points of a linear algebraic group, factoring out the solvable radical we may assume that $H$ is semi-simple and must prove that $H$ is compact. For algebraic groups over local fields of characteristic 0, being compact is equivalent to being anisotropic. So assume that $H$ is isotropic. Then by the Jacobson-Morozov theorem $H$ contains a copy either of $SL_2(K)$ or of $PGL_2(K)$, contradicting amenability of $H$. \eprsk


\noindent
{\bf Author's address:}\\
Institut de Math\'ematiques-Unimail\\
11 Rue Emile Argand\\
CH-2000 Neuch\^atel - SWITZERLAND\\
alain.valette@unine.ch

\end{document}